\documentclass[12pt,a4paper]{article}
\usepackage[utf8]{inputenc}
\usepackage[english]{babel}
\usepackage{amsmath}
\usepackage{amsfonts}
\usepackage{amssymb}
\usepackage{makeidx}
\usepackage{graphicx}
\usepackage[left=2cm,right=2cm,top=2cm,bottom=2cm]{geometry}
\usepackage{stmaryrd}
\usepackage{amsthm}
\usepackage{xcolor}
\usepackage[all]{xy}
\usepackage{tikz}

\providecommand{\keywords}[1]
{
  \small	
  \textbf{Keywords:} #1
}

\title{Cyclic splittings of pro-$\C$ groups}

\author{Jesus Berdugo $^{2}$, Pavel Zalesskii$^{1}$\footnote{Partially supporte by FAPDF and CNPq}  \\
        \small $^{1}$ Department of mathematics, University of Brasilia, E-mail: pz@mat.unb.br  \\
        \small $^{2}$ {Departamento de Matemáticas, Universidad del Atlántico,  E-mail: jeberdugo@mail.uniatlantico.edu.co} \\
}

\date{\today}

\newtheorem{teo}{Theorem}[section]

\newtheorem{theorem}[teo]{Theorem}
\newtheorem{proposition}[teo]{Proposition}%
\newtheorem{example}[teo]{Example}%
\newtheorem{remark}[teo]{Remark}%
\newtheorem{corollary}[teo]{Corollary}

\newtheorem{definition}[teo]{Definition}%
\newtheorem{lemma}[teo]{Lemma}

\newcommand{\Z}{\mathbb{Z}}
\newcommand{\G}{\mathcal{G}}
\newcommand{\C}{\mathcal{C}}
\begin{document}
\maketitle

\abstract{In this paper we proved a pro-$\C$ version of the Rips-Sela Theorems on  splittings as an amalgamated free product or HNN-extension over a cyclic subgroup, where $\C$ is a class of finite groups closed for subgroups, quotients, finite direct products and extensions with abelian kernel.

\vspace{0.5cm}


\section{Introduction}

In  1997 Rips and Sela published the fundamental paper \cite{bib1}, where they studied infinite cyclic splittings (i.e. $\Z$-splittings) of groups as an amalgamated free product or an HNN-extension.  They  constructed a canonical JSJ decomposition  for finitely presented groups  that gives a complete description of all $\Z$-splittings of these groups.
In order to understand all possible $\Z$-splittings of a group, they needed to study carefully the "interaction" between any two given elementary $\Z$-splittings of it.
 
The objective of this paper is to study the ``interaction" between any two given  $\widehat\Z_\C$-splittings of a pro-$\C$ group, where $\C$ is a class of finite groups closed for subgroups, quotient finite direct products and  extensions with abelian kernel. Namely, we prove a pro-$\C$ version of Rips-Sela's theorems on $\mathbb{Z}$-splittings ( \cite[Theorem 2.1 and Theorem 3.6]{bib1}).

 Let $\pi$ be a set of primes  and $\Z_\pi=\prod_{p\in\pi} \Z_p$ be a torsion free cyclic pro-$\pi$ group.
A splitting of a pro-$\C$ group $G$ as an amalgamated free pro-$\C$ product or HNN-extension over a cyclic group will be called a cyclic splitting and if this cyclic subgroup is torsion free, i.e.  isomorphic to  $\mathbb{Z}_\pi$ for some set of primes, then it will be called a $\mathbb{Z}_\pi$-splitting in the paper. An element $g$ of $G$ is called elliptic with respect to a cyclic splitting  as an amalgamated free pro-$\C$ product $G_1\amalg_C G_2$ (resp. pro-$\C$ HNN-extension $G=HNN(G_1, \Z_\pi, t)$) if $g$ is conjugate into $G_1\cup G_2$  (resp. into $G_1$) and is called hyperbolic if $\langle g \rangle$ intersects trivially any conjugate of $G_1$ and $G_2$ (resp. of $G_1$)  in $G$. A pair of given cyclic splittings $A_1\amalg_{C_1} B_1$ (resp. $G=HNN(G_1, C_1, t)$)) and $A_2\amalg_{C_2} B_2$ (resp. $G=HNN(G_2, C_2, t)$)) over $C_1=\langle c_1\rangle$, $C_2=\langle c_2\rangle$ is called:

\begin{itemize}

\item $Elliptic-Elliptic:$ If $c_{1}$ is elliptic in $A_2\amalg_{C_2} B_2$ and $c_{2}$ is elliptic in $A_1\amalg_{C_1} B_1$.

\item $Hyperbolic-Hyperbolic:$ If $c_{1}$ is hyperbolic in $A_2\amalg_{C_2} B_2$ and $c_{2}$ is hyperbolic in $A_1\amalg_{C_1} B_1$.

\item $Hyperbolic-Elliptic:$ If $c_{1}$ is hyperbolic in $A_2\amalg_{C_2} B_2$ and $c_{2}$ is elliptic in $A_1\amalg_{C_1} B_1$.
\end{itemize}

Note that in the abstract and pro-$p$ cases a non-elliptic element automatically means hyperbolic. In the pro-$\C$ situation (in particular, in  profinite) it is not always the case, but if the splitting is obtained from the pro-$\C$ completion of the abstract splitting then $C_i$ are either elliptic or hyperbolic (see Remark \ref{abstract hyperbolic is profinitely hyperbolic}). 
The Rips-Sela \cite[Theorem 2.1]{bib1} says that for a freely indecomposable group the Hyperbolic-Elliptic possibility does not exist. 

Our first result is the pro-$\C$ analog of it.

\begin{theorem}\label{sela} Let $G$ be a pro-$\C$ group that does not act on a pro-$\C$ tree with trivial edge stabilizers having no  global fixed point. Then $G$ does not admit a pair of hyperbolic-elliptic cyclic splittings. 
\end{theorem}

It is a well-known Bass-Serre's theory fact that a group splits as a free product if and only if it acts on  infinite tree  with trivial edge stabilizers (and without global fixed point). This is not true however in the pro-$\C$ case (in particular in the profinite case). It is known that the Kurosh Subgroup Theorem does not hold in the profinite case and in fact a non-cyclic subgroup $H$ of a free profinite product $G$ does not have to split at all as a non-trivial free profinite product (the most striking example of it is a Frobenius group $\Z_7\rtimes C_3$ \cite{GH}). However, $G$ acts on its standard profinite tree with trivial edge stabilizers and therefore so is $H$.  It follows that in the pro-$\C$  case the the hypothesis on $G$ in Theorem \ref{sela} is the pro-$\C$ version of the non-splitting into a free product.   

\begin{theorem}\label{selageneral} Let $G$ be a pro-$\C$ group that does not act on a pro-$\C$ tree with trivial edge stabilizers having no  global fixed point.  Then $G$ does not admit a pair of hyperbolic-elliptic cyclic splittings. 
\end{theorem}

To prove this theorem we  introduce two new constructions: the pro-$\C$ version  of blowing-up a vertex and the standard refinement (see Proposition \ref{refinement}).

The pro-$p$ version of Theorem \ref{sela} was proved in \cite{BZ} for finitely generated case. As Theorem \ref{selageneral} does not have this restricion we deduce the   pro-$p$ version of the Rips-Sela theorem without this restriction.

\begin{corollary} Let $G$ be a pro-$p$ group that does not act on a pro-$p$ tree with trivial edge stabilizers having no  global fixed points. Then any two $\mathbb{Z}_{p}$-splittings of $G$ are either elliptic-elliptic or hyperbolic-hyperbolic.

\end{corollary}

Next, given two hyperbolic-hyperbolic splittings over $C_1$ and $C_2$, we study the normalizer $N_G(C_i)$, $i=1,2$.  This study corresponds to the results of Section 3 of Rips-Sela's paper \cite{bib1} (where normalizers are called anti-centralizers). Note that Sela and Rips use  the existence of Tits' axis, on which the hyperbolic element  acts. In the pro-$\C$ case, such an axis does not exist, so our argument is different from their argument. 

\begin{proposition} \label{prop}  Let $G$ be pro-$\C$ group and  $G=A_i\amalg_{C_i} B_i$ (resp. $G=HNN(G_i, C_i, t)$) be a pair of  cyclic splittings. Let $H_i\neq 1$ be a subgroup of $C_i$. Then  $N_{G}(H_{i})$ is isomorphic  to one of the following groups:

\begin{enumerate}
\item[(i)]  $\Z_\pi\rtimes \Z_\rho$ ($\pi\cap \rho=\emptyset$);
\item[(ii)] $\mathbb{Z}_{\pi}\rtimes \Z/2$ infinite dihedral pro-$\C$ group;

\item[(iii)] a profinite Frobenius group $\Z_\pi\rtimes C_n$.

\item[(iv)] $ C_i\times \widehat\Z_\C$;

\item[(v)]   $N_{G}(H_i)=N_{A_i}(H_i)\amalg_{C_i}N_{B_i}(H_i)$, $[N_{A_i}(H_i):C_i]=2=[N_{B_i}(H_i):C_i]$ (this case can happen only  if $\Z/2\in \C$). Moreover, setting $\pi,\rho$ to be sets of primes such that $C_i\cong \Z_\pi$ and $\rho=\pi\setminus\{2\}$,  one of the following holds
\begin{enumerate}
\item[(a)] $N_{A_i}(H_i)$, $N_{B_i}(H_i)$ are torsion free and $N_G(H_i)$ is the pseudo Klein bottle  $(\Z_\rho\times \widehat\Z_\C)\rtimes \Z_2$ with $\Z_2$ acting by inversion, or 
\item[(b)] both $N_{A_i}(H_i)$ and $N_{B_i}(H_i)$ are  infinite dihedral  and
$N_G(H_i)$ is  isomorphic to $(\Z_{\pi}\times \widehat\Z_\C)\rtimes \Z/2$ with $\Z/2$ acting by inversion; 
\item[(c)] One of $N_{A_i}(H_i)$, $N_{B_i}(H_i)$ is torsion free and  the other  is  infinite dihedral, and   $N_G(H_i)\cong (\Z_\rho\times \widehat\Z_\C)\rtimes \Z_2)\rtimes \Z/2)$.
\end{enumerate}
In particular $N_G(H_i)$ contains an abelian normal subgroup of index $2$.
\end{enumerate}

\end{proposition}

If the class $\C$ consists of groups of odd order only then  the statement of the proposition is much simpler; we state it as the following

\begin{corollary} Suppose $\C$ consists of groups of odd order and $H_i\neq 1$ be a subgroup of $C_i$. Then  $N_{G}(H_{i})$ is isomorphic  to one of the following groups:
	
	\begin{enumerate}
		\item[(i)]  $\Z_\pi\rtimes \Z_\rho$ ($\pi\cap \rho=\emptyset$);

		\item[(ii)] a profinite Frobenius group $\Z_\pi\rtimes C_n$.

		\item[(iii)] $N_{G}(H_{i})\cong  C_i\times \widehat\Z_\C$.

		\end{enumerate}
	
	\end{corollary}

\bigskip
This description allows us to prove the second main result of the paper, namely the pro-$\C$ version of \cite[Theorem 3.6]{bib1}
 
\begin{theorem}\label{teo}
	Let $G$ be a  pro-$\C$ group  which can not  act on a pro-$\C$ tree with finite cyclic edge stabilizers  having no  global fixed point.  Let $G=A_{1}\amalg_{C_{1}} B_{1}$ (or $G=HNN(A_{1},C_{1},t_{1})$), and $G=A_{2}\amalg_{C_{2}} B_{2}$ (or $G=HNN(A_{2},C_{2},t_{2})$) be two hyperbolic-hyperbolic cyclic splittings of $G$. Suppose that $N_{G}(C_1)$ is not virtually cyclic. Then  $G=N_G(C_1)=N_G(C_2)$ is virtually abelian (isomorphic to one of  pro-$\C$ groups listed in (iv) or (v) of Proposition 1.4 ).  Moreover, if $\C$ consists of groups of odd order, then $G\cong C_i\times \widehat\Z_\C$.
	
\end{theorem} 
 
 The theorem shows in particular that if $N_G(C_1)$ is not virtually cyclic and G does not split (as an amalgamation or HNN-extension) over a finite cyclic subgroup, then $N_G(C_1)$ is in fact the ambient group $G$; hence is virtually  abelian.
 
 As a corollary we deduce an infinitely generated pro-$p$ version of the Rips-Sela result generalizing \cite[Theorem 5.2 and Theorem 5.11]{BZ}.
 
 \begin{corollary} Let $G$ be a  pro-$p$ group  which can  not act on a pro-$p$ tree with  edge stabilizers of order $\leq 2$ having no  global fixed point. Let $G=A_{1}\amalg_{C_{1}} B_{1}$ (or $G=HNN(A_{1},C_{1},t_{1})$), and $G=A_{2}\amalg_{C_{2}} B_{2}$ (or $G=HNN(A_{2},C_{2},t_{2})$) be two hyperbolic-hyperbolic $\mathbb{Z}_{p}$-splittings of $G$. Suppose that $N_{G}(C_1)$ is not  cyclic or infinite dihedral. Then  $G=N_G(C_1)=N_G(C_2)$ is  abelian of rank 2 by at most cyclic of order $2$. In particular, if $p>2$ then $G$ is abelian. 
 
 \end{corollary}

 The proofs of the results use  the pro-$\C$ version of the Bass-Serre theory  that can be found in \cite{bib3}.\\ 
 
The structure of the paper is as follows: in Section \ref{sec2} we give basic definitions that will be used throughout the article and introduce the pro-$\C$ version of blowing-up a vertex and standard refinement. In Section \ref{sec4} we prove Theorem \ref{sela} and  in  Section \ref{sec6} we prove Theorem \ref{teo}. In Section 6 we show that two hyperbolic-hyperbolic splittings such that $C_1$ does not commute with a conjugate of $C_2$ act 2-acylindrically on a pro-$\C$ tree.
 
 \bigskip
 {\bf Conventions.} Throughout the paper, unless otherwise stated, groups are profinite, subgroups are closed, and homomorphisms are continuous. In particular, $\langle S\rangle$ will mean the topological generation in the paper. $\C$ stands for a class of finite groups closed for subgroups, quotients, finite direct products and extensions with abelian kernel.  The presentations are taking in the category of pro-$\C$ groups; $a^g$ will stand for $g^{-1}ag$ in the paper. $\widehat\Z_\C$ will denote the pro-$\C$ completion of $\Z$ and $\Z_\pi$ the pro-$\pi$ completion (i.e. $\Z_\pi=\prod_{p\in\pi} \Z_p$) of $\Z$, where $\pi$ is a set of primes.

\section{Preliminaries}\label{sec2}
\subsection{Basic definitions (pro-$\C$ Bass-Serre' theory)}

Unless otherwise specified, in this section $\C$ is a class of finite groups closed for subgroups, quotients, finite direct products and extensions with abelian kernel. In this section we will give necessary definitions and known results of the pro-$\C$ version of the Bass-Serre theory.
 
\begin{definition}
A $graph$ $\Gamma$ is a disjoint union $E(\Gamma)\cup V(\Gamma)$, with two maps $$d_0,d_1:\Gamma\longrightarrow V(\Gamma),$$ whose restrictions to $V(\Gamma)$ are the identity map.   For any element $e\in E(\Gamma)$, $d_0(e)$ and $d_1(e)$ are called the initial and the terminal vertices of $e$ respectively. 
\end{definition}

\begin{definition}\label{prograph}
$\Gamma$ is called a profinite graph  if $\Gamma$ is a profinite space such that:

\begin{itemize}
    \item[1.] $V(\Gamma)$ is closed,
    \item[2.] the maps $d_0,d_1:\Gamma\longrightarrow V(\Gamma)$ are continuous.
\end{itemize}
\end{definition}

We call
$V (\Gamma) := d_0(\Gamma) \cup d_1(\Gamma)$ the set of vertices of $\Gamma$ and $E(\Gamma) := \Gamma \setminus V (\Gamma)$ the set of edges of $\Gamma$. For $e \in E(\Gamma)$ we call $d_0(e)$ and $d_1(e)$  the initial and
terminal vertices of  the edge $e$. 

A \emph{morphism} $\alpha:\Gamma\longrightarrow \Delta$ of profinite graphs is a continuous map with $\alpha d_i=d_i\alpha$ for $i=0,1$. If $\alpha$ is injective, the image  is called a subgraph of profinite graph $\Delta$, if $\alpha$ is surjective, then $\Delta$ is called a quotient graph of $\Gamma$. Note that a morphism of graphs can send edges to vertices.

A profinite graph is called connected if every finite quotient graph of it is connected. A maximal connected profinite  subgraph is called a connected component.

\begin{definition}\label{collapse}
	 Let $\Delta$  be a profinite subgraph of a profinite graph $\Gamma$. 
	Define the operation of
collapsing the connected profinite components of $\Delta$
to points as 
a natural mapping
to the quotient space 
$\Gamma/\sim$, 
where $\sim$ is the equivalence relation defined as
follows: $m\sim m'$ if  either \( m = m' \), or \( m, m' \in \Delta^*(m) \), where $\Delta^*(m)$ is the connected component of $\Delta$ containing $m$. 
	Then \( \alpha \) is a morphism of graphs, and \( \Gamma / \Delta \) becomes a quotient graph of \( \Gamma \). 
	We say that \( \Gamma / \Delta \) is obtained from \( \Gamma \) by \emph{collapsing} the connected components of \( \Delta \) to a point.  
	Then the collapse of all connected components of 
	$\Delta$ is an epimorphism of profinite 
	graphs 
	(see \cite[Exercise 2.1.11]{bib3})
	.

\end{definition}

 A profinite graph $\Gamma$ is said to be a pro-$\C$ tree, if the following sequence of profinite $\widehat{\mathbb{Z}}_\C$-modules is exact

$$ 0 \longrightarrow \llbracket\widehat\Z_\C(E^{\ast}(\Gamma),\ast) \rrbracket {\buildrel\delta\over\longrightarrow} \llbracket \widehat\Z_\C( V(\Gamma)) \rrbracket{\buildrel\epsilon\over\longrightarrow} \widehat\Z_\C \longrightarrow 0,$$
where  $E^{\ast}(\Gamma)=\Gamma/V(\Gamma)$, and the homomorphisms are defined $\delta(\bar{e})=d_1(e)-d_0(e)$, with $\bar{e}$ be the image of edge $e\in E(\Gamma)$ in $E^{\ast}(\Gamma)$, and $\delta(\ast)=0$; and $\epsilon(v)=1$, for all $v\in V(\Gamma)$.\\

Now consider $\pi=\pi(\mathcal{C})$ to be the set of prime numbers involved in the pseudovariety $\mathcal{C}$, i.e, the set of all prime numbers $p$ that divide the order of some group $H$ in $\mathcal{C}$. The notion of a pro-$\C$ tree depends only on $\pi$ so we can speak about pro-$\pi$ trees  if it is more convenient.

If $v,w\in V(\Gamma)$, the smallest pro-$\pi$ subtree of $\Gamma$ containing $\{v,w\}$ is called the $geodesic$ connecting $v$ and $w$, and is denoted $[v,w]$ (the definition is in pag. 56 of \cite{bib3}). By definition a pro-$\C$ group $G$ $acts$ on a profinite graph $\Gamma$ if we have a continuous action of $G$ on the profinite space $\Gamma$, such that $d_0$ and $d_1$ are $G$-maps. We say that the action of a pro-$\C$ group $G$ on a pro-$C$ tree $T$ is $irreducible$, if $T$ is the unique minimal $G$-invariant pro-$\C$ subtree of $T$. The action is said to be $faithful$ if the kernel $K$ of the action is trivial. If $G$ acts on $T$ irreducibly then the resulting action of $G/K$ on $T$ is faithful and irreducible.

\begin{proposition}\label{minimal subtree}(\cite[Lemma 1.5]{Zal90} or \cite[Proposition 2.4.12]{bib3}).  Let $G$ be a pro-$\C$ group acting on a pro-$\C$ tree $T$. Then there exists a non-empty minimal $G$-invariant subtree $D$ of $T$. Moreover, if $G$ does not stabilize a vertex, then $D$ is unique.

\end{proposition}

\begin{definition} \rm
 Let $G$ be a pro-$C$ group, and $T$  a pro-$\C$ tree  on which $G$ acts continuously. For $g\in G$,
 \begin{itemize}
 \item $g$ is $elliptic$ if it fixes a vertex in $T$; if a subgroup of $G$ fixes a vertex we also call it elliptic.
 \item $g$ is $hyperbolic$ if the subgroup $\langle g\rangle$ generated by $g$ acts freely on $T$. 
 \end{itemize}
 \end{definition}

 If $g$ is hyperbolic, then by  Proposition \ref{minimal subtree} there exist a unique nonempty minimal $\langle g\rangle$-invariant pro-$\pi$ subtree $D_g\subseteq T$.

\begin{remark}\label{abstract hyperbolic is profinitely hyperbolic}In the profinite case if an element $g$ is not elliptic it is not necessarily hyperbolic, because $\langle g\rangle$ might not act freely on $T$. However, if the action of $G$ on a profinite tree comes from the action of the abstract group $G^{abs}$  on a tree $T^{abs}$ and if $b\in G^{abs}$ is hyperbolic in $T^{abs}$ it is hyperbolic in $T$ (see \cite[Lemma 8.3.2]{R} for details). 
\end{remark}

  \begin{example} Let $G=G_1*_C G_2$ be a free product of residually $\C$-groups with infinite cyclic amalgamation $C$. Then the pro-$\C$ completion $\widehat G=\widehat G_1\amalg_{\widehat C} \widehat G_2$ is a pro-$\C$ splitting of $\widehat G$ over $\widehat C$ (see \cite{R}). If $b\in G$, by  \cite[Proposition 2.9]{RZ-96} or \cite[Lemma 8.3.2]{bib3} $b$ is hyperbolic in $G$ if and only if it is hyperbolic in $\widehat G$. This justifies the definition above. 
\end{example}

Note that in the pro-$p$ case an element is hyperbolic if and only if it is not elliptic. In the pro-$\C$ case we can have a mixed case; some of the $p$-components of $g$ are elliptic and some of them are hyperbolic. For example, $C_2\times \Z_3$ can act on a profinite tree with $C_2$ acting trivially and $\Z_3$ acting freely. Then the generator $(a,b)$ of $C_2\times \Z_3$, $1\neq a\in C_2$, $b\in \Z_3$ is an element of  mixed type.

\medskip 
 The following lemma is used to  prove that the graph we construct in Proposition \ref{refinement} is a pro-$\C$ tree.

\begin{lemma}\label{Lemma.standard.help}  Let $T$ be a connected profinite graph and $D$ is a subgraph of $T$ such that each connected component of it is a pro-$\mathcal{C}$ tree. Let $\widehat{T}$ be the profinite quotient graph of  $T$ obtained by collapsing all connected components of $D$. Then  $\widehat{T}$ is a pro-$\mathcal{C}$ tree if and only if $T$ is a pro-$\mathcal{C}$ tree.
\end{lemma}

\begin{proof} Let 
$$ 0 \longrightarrow \llbracket(E^{\ast}(T),\ast) \rrbracket {\buildrel\delta\over\longrightarrow} \llbracket \widehat\Z_\C V(T) \rrbracket{\buildrel\epsilon\over\longrightarrow} \widehat\Z_\C \longrightarrow 0,$$
and 

$$ 0 \longrightarrow \llbracket\widehat\Z_\C(E^{\ast}(\widehat{T}),\ast) \rrbracket {\buildrel\hat\delta\over\longrightarrow} \llbracket \widehat\Z_\C V(\widehat{T}) \rrbracket{\buildrel\hat\epsilon\over\longrightarrow} \widehat\Z_\C \longrightarrow 0,$$
be the sequences associated to $T$ and $\widehat T$. Note that the coonectedness of $T$ and $\hat T$ is equivalent to the exactness of these sequences in the middle term  by \cite[Proposition 2.3.2]{bib3}.
	Consider $\Delta=\{D^{\ast}(m):m\in D\}$, the set of all connected component of $D$. Since $D^{\ast}(m)$ is a pro-$\mathcal{C}$ tree for each $m\in D$, the map
	 $$0\longrightarrow \llbracket \widehat\Z_\C(E^{\ast}(D^{\ast}(m)),\ast) \rrbracket\xrightarrow{\delta} \llbracket \widehat\Z_\C V(D^{\ast}(m)) \rrbracket$$ is exact for all $m\in D$. Hence  $0\longrightarrow \llbracket \widehat\Z_\C(E^{\ast}(D,\ast) \rrbracket \xrightarrow{\delta_{\varDelta}} \llbracket \widehat\Z_\C V(D) \rrbracket$ is exact (because $Ker(\delta_\Delta)\leq Ker(\delta)$). 
	 
	 Suppose  $\widehat{T}$ is a pro-$\C$ tree, so that the bottom sequence is exact. Let  $\beta:D \longrightarrow T$ be the natural inclusion  and $\alpha: T\longrightarrow \widehat{T}$ the natural epimorphism. Then we have the following commutative diagram: 
	$$
	\xymatrix{ 0 \ar[r] & \llbracket \widehat\Z_\C(E^{\ast}(D,\ast) \rrbracket \ar[r]^{{\delta_{\varDelta}}} \ar[d]^{\tilde{\beta}_{E^{\ast}}} &  \llbracket \widehat\Z_\C V(D) \rrbracket \ar[d]^{\tilde{\beta}_{V}} & & \\
		& \llbracket \widehat\Z_\C(E^{\ast}(T),\ast) \rrbracket \ar[r]^{\delta} \ar[d]^{\tilde{\alpha}_{E^{\ast}}} & \llbracket \widehat\Z_\C V(T) \rrbracket \ar[r]^{\epsilon} \ar[d]^{\tilde{\alpha}_{V}} & \widehat\Z_\C\ar[r] \ar[d]^{Id_{\widehat\Z_\C}} & 0\\
		0\ar[r] & \llbracket \widehat\Z_\C(E^{\ast}(\widehat{T}),\ast) \rrbracket \ar[r]^{\hat\delta} & \llbracket \widehat\Z_\C V(\widehat{T}) \rrbracket \ar[r]^{\hat\epsilon} & \widehat\Z_\C \ar[r]& 0  }
	$$ \\
	Clearly $(E^{\ast}(\widehat{T}),\ast)\cong (E^{\ast}(T),\ast)/(E^{\ast}(D),\ast) $, thus by \cite[Lemma 1.9]{bib6} $$\llbracket \widehat\Z_\C(E^{\ast}(\widehat{T}),\ast) \rrbracket\cong \llbracket \widehat\Z_\C(E^{\ast}(T),\ast) \rrbracket/ \llbracket \widehat\Z_\C(E^{\ast}(\Delta),\ast) \rrbracket.$$ 
	Then  $Ker(\tilde{\alpha}_{E^{\ast}})=\tilde{\beta}_{E^{\ast}}(\llbracket \widehat\Z_\C(E^{\ast}(D),\ast) \rrbracket)$. Take $a\in Ker(\delta)$. Since $\hat\delta\tilde{\alpha}_{E^{\ast}}(a)=\tilde{\alpha}_{V}\delta(a)=0$,  $\tilde{\alpha}_{E^{\ast}}(a)\in Ker(\hat\delta)=0$ and so $a\in Ker(\tilde{\alpha}_{E^{\ast}})=\tilde{\beta}_{E^{\ast}}(\llbracket \widehat\Z_\C(E^{\ast}(D),\ast) \rrbracket)$, i.e, there exists $b\in \llbracket \widehat\Z_\C(E^{\ast}(D),\ast) \rrbracket$, such that $\tilde{\beta}_{E^{\ast}}(b)=a$. By commutativity of the diagram at the top, $\tilde{\beta}_{V}({\delta_{\varDelta}}(b))=\delta(\tilde{\beta}_{E^{\ast}}(b))=\delta(a)=0$. Since  $\tilde{\beta}_{V}$ and ${\delta_{\varDelta}}$ are injections, we have that  $b=0$, therefore $a=0$. Thus $Ker(\delta)=0$, i.e, $\delta$ is an injection. Hence $T$ is a pro-$\mathcal{C}$ tree ( since $Im(\delta)=Ker(\epsilon)$, $T$ is connected profinite graph by \cite[Proposition 2.3.2]{bib3}).

	Now assume that $T$ is a pro-$\mathcal{C}$ tree, so that the middle sequence is exact. Since $\widehat{T}$ is a quotient of $T$, then it is connected, hence $Im(\hat{\delta})=Ker(\hat{\epsilon)}$ by \cite[Proposition 2.3.2]{bib3}. Thus we have the following commutative diagram:

		$$
	\xymatrix{ 0 \ar[r] & \llbracket \widehat\Z_\C(E^{\ast}(D,\ast) \rrbracket \ar[r]^{\delta_{\varDelta}} \ar[d]^{\tilde{\beta}_{E^{\ast}}} &  \llbracket \widehat\Z_\C V(D) \rrbracket \ar[d]^{\tilde{\beta}_{V}} & & \\
	0\ar[r]	& \llbracket \widehat\Z_\C(E^{\ast}(T),\ast) \rrbracket \ar[r]^{\delta} \ar[d]^{\tilde{\alpha}_{E^{\ast}}} & \llbracket \widehat\Z_\C V(T) \rrbracket \ar[r]^{\epsilon} \ar[d]^{\tilde{\alpha}_{V}} & \widehat\Z_\C\ar[r] \ar[d]^{Id_{\widehat\Z_\C}} & 0\\
		 & \llbracket \widehat\Z_\C(E^{\ast}(\widehat{T}),\ast) \rrbracket \ar[r]^{\hat\delta} & \llbracket \widehat\Z_\C V(\widehat{T}) \rrbracket \ar[r]^{\hat\epsilon} & \widehat\Z_\C \ar[r]& 0  }
	$$ \\
		Observe that $Ker(\widetilde{\alpha}_{V})=Ker(\epsilon)\cap \tilde{\beta}_{V}( \llbracket \widehat\Z_\C V(D) \rrbracket )=\delta(\llbracket \widehat\Z_\C(E^{\ast}(T),\ast) \rrbracket)\cap \tilde{\beta}_{V}(\llbracket \widehat\Z_\C V(D) \rrbracket) = \tilde{\beta}_{V}(\delta_{\varDelta}(\llbracket \widehat\Z_\C(E^{\ast}(D,\ast) \rrbracket)) .$
		
		
Choose $b\in Ker({\hat{\delta}})$. There exists $a\in  \llbracket \widehat\Z_\C(E^{\ast}(T,\ast) \rrbracket$ such that $\widetilde{\alpha}_{E^{\ast}}(a)=b$. Hence $\delta(a)\in Ker(\widetilde{\alpha}_{V})=\tilde{\beta}_{V}({\delta_{\varDelta}}(\llbracket \widehat\Z_\C(E^{\ast}(D,\ast) \rrbracket))$. Therefore there exists  $c\in \llbracket \widehat\Z_\C(E^{\ast}(D,\ast) \rrbracket $ such that $\tilde{\beta}_{V}({\delta_{\varDelta}}(c))=\delta(a)$, i.e, $\delta(\widetilde{\beta}_{E^{\ast}}(c))=\delta(a)$. Then $\widetilde{\beta}_{E^{\ast}}(c)-a\in Ker(\delta)=0$ and so $\widetilde{\beta}_{E^{\ast}}(c)=a$. Moreover $b=\widetilde{\alpha}_{E^{\ast}}(a)=\widetilde{\alpha}_{E^{\ast}}(\widetilde{\beta}_{E^{\ast}}(c))=0$, so $b=0$, thus $Ker(\hat\delta)=0$.
\end{proof}

\bigskip

\begin{definition}
Let $\Delta$ be a connected finite graph. A \textit{graph of pro-$\C$ groups} $(\mathcal{G},\Delta)$ over $\Delta$ consists of a pro-$\C$ group $\mathcal{G}(m)$ for each $m\in \Delta$ and continuous monomorphisms $\partial_i: \mathcal{G}(e)\longrightarrow \mathcal{G}(d_i(e))$ for each edge $e\in E(\Delta)$, $i=0,1$.
\end{definition}

\begin{definition} The   fundamental group $\Pi^{\C}_{1}(\mathcal{G},\Delta)$ of a finite graph of pro-$\C$ groups $(\mathcal{G},\Delta)$ is a completion of the abstract fundamental group $\pi^{abs}_{1}(\mathcal{G},\Delta)$ of this graph of groups $(\mathcal{G},\Delta)$ (for the definition of $\pi^{abs}_{1}(\mathcal{G},\Delta)$  see \cite{bib2}), in the following way: consider the collection $\mathcal{J}=\{J\trianglelefteq_{f}\pi_{1}^{abs}(\mathcal{G},\Delta):J\cap \mathcal{G}(v)\trianglelefteq_{o}\mathcal{G}(v)\} $ then 

 $$\Pi^{\C}_{1}(\mathcal{G},\Delta)=\lim_{\overleftarrow{J \in\mathcal{J}}}\pi^{abs}_{1}(\mathcal{G},\Delta)/J. $$  
 \end{definition}
 
Note that in contrast with the abstrtact case, the groups $\G(m)$ not always embed in $\Pi^{\C}_{1}(\mathcal{G},\Delta)$, but one can always replace them by their images in $\Pi^{\C}_{1}(\mathcal{G},\Delta)$ and then the embedding will hold. The fundamental group  $\Pi^{\C}_{1}(\mathcal{G},\Delta)$ will not change by this operation. Thus from now on we shall always assume that such embeddings hold, i.e. the graph of groups is injective in the terminology of \cite{bib3}.
 
 \begin{example} If $\Delta$ consists of one edge $e$ and two vertices $v,w$ then $\Pi^{\C}_{1}(\mathcal{G},\Delta)=\G(v)\amalg_{\G(e)}\G(w)$ is an amalgamated free pro-$\C$ product. Note that  the factors do not always embed in an amalgamated free pro-$\C$ products,  i.e. an amalgamated free pro-$\C$ product is not always  proper (in the terminology of  \cite{bib4}). However, we can make it proper by replacing $\G(v),\G(w),\G(e)$ with their images in $G$ (as explained in the subsection 9.2 and 9.4 in \cite{bib5}).
 
 \end{example}
 
 \begin{example} If $\Delta$ consists of one edge $e$ and one vertex $v$ then $$\Pi^{\C}_{1}(\mathcal{G},\Delta)=HNN(\G(v), \G(e), t)$$ is a pro-$\C$ HNN-extension.  Note that   the base group does not always embed in the pro-$\C$ HNN-extension,  i.e.  the pro-$\C$ HNN-extension is not always  proper (in the terminology of  \cite{bib4}). However, we can make it proper by replacing $\G(v),\G(e)$ with their images in $G$ (as explained in the subsection 9.2 and 9.4 in \cite{bib4}).
 
 \end{example}

\begin{definition}
The standard pro-\( C \) tree (or universal covering graph) of a finite graph of pro-\( C \) groups \( (\mathcal{G}, \Delta) \) is defined as:
\[
S = S(\mathcal{G}) = \bigcup_{m \in \Delta} G / \mathcal{G}(m)
\]
(cf. \cite[Example 3.6.1]{bib3} and \cite[Theorem 6.5.2]{bib3}). The vertices of \( S \) are those cosets of the form \( g\mathcal{G}(v) \), with \( v \in V(\Delta) \) and \( g \in G=\Pi^{\C}_{1}(\mathcal{G},\Delta) \); its edges are the cosets of the form \( g\mathcal{G}(e) \), with \( e \in E(\Delta) \); and the incidence maps of \( S \) are given by the formulas
\[
d_0(g\mathcal{G}(e)) = g\mathcal{G}(d_0(e)); \quad 
d_1(g\mathcal{G}(e)) = g t_e \mathcal{G}(d_1(e)) \quad 
(e \in E(\Delta)).
\]     
\end{definition}

The fundamental group $\Pi^{\C}_{1}(\mathcal{G},\Delta)$ of a finite graph of pro-$\C$ groups $(\mathcal{G},\Delta)$ acts on a standard pro-$\C$ tree.

We write explicitly the standard trees $S(G)$ on which $G$ acts for the cases of an amalgamated free pro-$\C$ product $G=G_1\amalg_HG_2$ 
  and an HNN-extension $G=HNN(G_1, H,t)$  since we shall use them very often for these cases.

  \begin{itemize}
  \item Let $G=G_1\amalg_HG_2$. Then
 the vertex set is $V(S(G))= \displaystyle G/G_1\cup G/G_2$,
  the edge set is $E(S(G))= G/H$, and
  the initial and terminal vertices of an edge $gH$ are
  respectively  $gG_1$ and $gG_2$.
   
    \item Let $G=HNN(G_1, H,t)$. Then
 the vertex set is $V(S(G))= \displaystyle G/G_1$,
  the edge set is $E(S(G))= G/H$, and
  the initial and terminal vertices of an edge $gH$ are
  respectively  $gG_1$ and $gtG_1$.
  \end{itemize}

We shall often use the following 

\begin{proposition}\label{intersection}(\cite[Corollary 7.1.5]{bib3} or \cite[Theorem 3.12]{ZM})
\begin{enumerate}
\item[(a)]  Let $G=G_1\amalg_H G_2$ be a proper free pro-$\C$ product with amalgamation. Then $$G_i\cap G_j^g\leq H^b$$ for some $b\in G_i$, whenever $i\neq j$ or $g\not\in G_i$ ($i,j=1,2$).

\item[(b)] Let $G=HNN(G_0, H, t)$ be a proper pro-$\C$ HNN-extension. Then  $$G_0\cap (G_0)^{g^{-1}}\leq A^b\  {\rm or}\  G_0\cap (G_0)^{g^{-1}}\leq (A^t)^b$$ for some $b\in G_0$.
\end{enumerate}

\end{proposition}

We finish the subsection  presenting the {\em blowup} of a subset $V$ of vertices of a profinite graph $\Gamma$. 
The idea consists of replacing each of these vertices by the connected components
of a given profinite  graph $\tilde\Sigma$ -- yielding a profinite subgraph inside the blowup. Since collapsing each of its components gives us back $\Gamma$ we consider
the construction as a reverse procedure to {\em collapsing the connected components of a  subgraph}. 
We describe the construction in a more formal way.

\begin{definition} Let $G$ be a profinite group and $H\leq G$ be a subgroup of $G$. If $H$ acts on a profinite graph $\Gamma$, the induced $G$-graph $G\times_H \Gamma$ is defined to be $G\times \Gamma$ (with $G$ acting on the first coordinate by multiplication)  modulo the equivalence relation $(gh,m)\sim (g,hm)$ $g\in G, h\in H, m\in \Gamma$.
\end{definition}

\begin{definition} \label{standard refinement}
	
	Let  $(\G,\Gamma)$ be a finite graph of pro-$\C$ groups and $v\in V(\Gamma)$ be a vertex. Suppose $G(v)$ acts on a pro-$\C$ tree $T$ irreducibly   such that the edge group $G(e)$ of every edge incident to $v$ stabilizers a vertex $v_e$ in $T$. Let $S(G)$ be the standard pro-$\C$ tree of $G=\Pi_1^{\C}(\G,\Gamma)$. One replaces the set of vertices $G/G(v)\subseteq  S(G)$ by $G/G(v)\times_{G(v)} T$ and one attaches for each  edge $e$ incident to $v$ the edge $gG(e)$  to the  vertex $(g, v_e)\in G/G(v)\times_{G(v)} T$, i.e. $d_i(gG(e))=(g, v_e)$ if  $d_i(e)=v$.  The obtained  pro-$\C$ graph $\widehat T$ is called a standard refinement of $(\G,\Gamma)$ or of $S(G)$ obtained by blowing up the vertex $v$.\end{definition}

\begin{proposition}\label{refinement} Let   $S(G)$  be  standard graph for $(\G,\Gamma)$ and \(v\) a vertex in \(\Gamma\).  Suppose that $G(v)$ acts on a pro-$\C$ tree $T$ irreducibly   such that the edge group $G(e)$ of every edge incident to $v$ stabilizers a vertex $v_e$ in $T$. Then there exists a profinite graph $\widehat{S}$ and a morphism $p:\widehat S \longrightarrow S(G)$  such that
	
	\begin{enumerate}
		\item[(1)] $p$ is a collapse map (cf. Definition \ref{collapse});
		
		\item[(2)] every vertex stabilizer of $\widehat S$  fixes a vertex  either in  $T$ or in $S(G)$;
		\item[(3)] the stabilizer of any edge of $\widehat S$ fixes an edge in $S(G)$ or in $T$;
	\end{enumerate}		
	In particular, $\widehat S$ is a pro-$\C$ tree.
\end{proposition}

\begin{proof}  
	We define $p:\widehat{S}\longrightarrow S(G)$ in standard manner that sends  $gT$ to $gv$. This map clearly satisfy the first requirement.
	
	Let us check that the other properties follow. Assertion (2) is clear (by construction of $\widehat{T}$).
	
	If $e$ is an edge of $\widehat{S}$ that is not collapsed by $p$, then $G_e$ fixes an edge $p(e)$ of $S(G)$. Otherwise $e\in T$ so $G_e$ fixes an edge in $T$, and Assertion (3) holds.
	
	The last statement follows from Lemma \ref{Lemma.standard.help}.

\end{proof}

\subsection{ Cyclic splittings}

 \bigskip
Let $G$ be a pro-$\C$ group. A cyclic splitting of $G$ is a splitting as non-fictitious proper free pro-$\C$ product with infinite cyclic amalgamation or as a proper pro-$\C$ HNN-extension with infinite cyclic associated subgroup. If this cyclic subgroup is torsion free, i.e.  isomorphic to $\Z_\pi$ for some set of primes $\pi$, we call the splitting $\Z_\pi$-splitting. Let $C_1=\langle c_1\rangle $ and $C_2=\langle c_2\rangle$ be cyclic subgroups of $G$  and suppose that $G$ admits   two cyclic splittings over them: 
\begin{itemize}
\item $G=A_{1}\amalg_{C_{1}} B_{1}$ or $G=HNN(A_{1},C_{1},t_{1}))$.
\item $G=A_{2}\amalg_{C_{2}} B_{2}$ or $G=HNN(A_{2},C_{2},t_{2}))$.
\end{itemize}
 
Let $T_1$ be the standard pro-$\C$ tree  corresponding to the first splitting and $T_2$ the standard pro-$\C$ tree corresponding to the second splitting. The pair of these  splittings is called:

\begin{itemize}

\item $Elliptic-Elliptic:$ If $c_{1}$ is elliptic in $T_{2}$ and $c_{2}$ is elliptic in $T_{1}$.

\item $Hyperbolic-hyperbolic:$ If $c_{1}$ is hyperbolic in $T_{2}$ and $c_{2}$ is hyperbolic in $T_{1}$.

\item $Hyperbolic-elliptic:$ If $c_{1}$ is hyperbolic in $T_{2}$ and $c_{2}$ is elliptic in $T_{1}$.
\end{itemize}

\subsection{Normalizer}
Item (1) of the following Proposition has been proved in \cite[Proposition 15.2.4]{bib3} in the profinite case and both items proved in  \cite[ Proposition 8.1]{bib7} for the pro-$p$ case. All these proofs go mutadis mutandis in the pro-$\C$ case, where $\C$ is a class of finite groups closed for subgroups, quotients, direct products and extensions with abelian kernel.

\begin{proposition} \label{teo6} Let $H$ be a non-trivial subgroup of a  cyclic subgroup 
	$C\cong\Z_\pi$. Then:

\begin{itemize}

\item[(1)] If $G=A\amalg_{C} B$ is a cyclic splitting as an amalgamated free pro-$\C$ product, then $N_{G}(H)=N_{A}(H)\amalg_{C} N_{B}(H)$.

\item[(2)] If $G=HNN(A,C,t)$ is a cyclic splitting as an HNN-extension then:

\begin{itemize}

\item If $C$ and $C^{t}$ are conjugate in $A$ then: $N_{G}(H)=HNN(N_{A}(H),C,t')$ and $G=HNN(G,C,t')$.

\item If $C$ and $C^{t}$ are not conjugate in $A$ then $N_{G}(H)=N_1\amalg_{C} N_{2}$, where $N_{1}=N_{A^{t^{-1}}}(H)$ and $N_{2}=N_{A}(H)$.

\end{itemize}

\end{itemize}

\end{proposition}

\begin{corollary}\label{centralizer} In addition to the hypotheses of Proposition \ref{teo6}, suppose that $N_A(C)\neq N_G(C)\neq N_B(C)$ in the case $G=A\amalg_C B$. Then
	$C_G(C)$ contains $F\times C$, where $F$ is a free pro-$\C$ group. 
	
\end{corollary}

\begin{proof} 
	
	Case 1.  $N_G(C)=N_1\amalg_C N_2$, where $N_1=N_A(C)$ and $N_2=N_B(C)$ (if $G=A\amalg_C B$) or $N_1=N_{A^{t^{-1}}}(C)$ and $N_{2}=N_{A}(C)$  (if $G=HNN(A,C,t)$). The centralizer $C_G(C)$ is the kernel of the natural homomorphism $f:N_G(C)\longrightarrow Aut(C)$. Clearly, it contains $C$. It suffices to show that $C_G(C)/C$ contains a free pro-$\C$ group $F$. Note that $f$ factors through $N_G(C)/C=N_1/C\amalg N_2/C$ so it suffices to show that the kernel of the corresponding homomorphism $\bar f:N_G(C)/C\longrightarrow Aut(C)$ contains $F$.  But $Aut(C)$ is abelian, and so $\bar f$ factors through $N_1/C\times N_2/C$. The kernel of  the natural epimophism  $$N_G(C)/C=N_1/C\amalg N_2/C \longrightarrow N_1/C\times N_2/C$$ is the Cartesian subgroup which is free pro-$\C$ group $F$ on the pointed profinite space  $\{[n_1,n_2], n_1\in N_1/C, n_2\in N_2/C\}$ (see \cite[Theorem 9.1.6]{bib4} ).  The preimage of $F$ in $N_G(C)$ is isomorphic to $F\times C$ as required.

	Case 2.  $N_{G}(C)=HNN(N_{A}(C),C,t')$ and $G=HNN(G,C,t')$. Again the natural homomorphism $f:N_G(C)\longrightarrow Aut(C)$   factors through $N_G(C)/C=N_A(C)/C\amalg \widehat\Z_\C$ so it suffices to show that the kernel of the corresponding homomorphism $$\bar f:N_G(C)/C\longrightarrow Aut(C)$$ contains $F$.  But $Aut(C)$ is abelian, and so $\bar f$ factors through $N_A(C)/C\times \widehat\Z_\C$. The kernel of  the natural epimophism of $$N_G(C)/C=N_A(C)/C\amalg \widehat\Z_\C\longrightarrow N_A(C)/C\times \widehat\Z_\C$$ is the Cartesian subgroup which is free pro-$\C$ group $F$ on the pointed profinite space  $\{[n_1,n_2], n_1\in N_A(C)/C, n_2\in \widehat\Z_\C\}$ (see \cite[Theorem 9.1.6]{bib4} ).  The preimage of $F$ in $N_G(C)$ is isomorphic to $F\times C$ as required.
\end{proof}

\begin{corollary}\label{soluble normalizer}
 Let $H$ be a non-trivial subgroup of a  cyclic subgroup	$C\cong\Z_\pi$. Suppose $N_G(H)$ is soluble. Then either $C=N_A(H)$  or $C=N_B(H)$  or $[N_A(H):C]=2=[N_B(H):C]$. Moreover, if $G=HNN(A,C,t)$ then $N_G(H)=C\rtimes \widehat\Z_C$.

\end{corollary}

\begin{proof} Note first that since $C$ is cyclic,  $N_G(C)\leq N_G(H)$. A non-fictitious free amalgamated pro-$\C$ product $N_A(H)\amalg_C N_B(H)$ will not have non-abelian free pro-$\C$ subgroup, and so will not be soluble, only if $[N_A(H):C]=2=[N_B(H):C]$. 

Suppose now $G=HNN(A,C,t)$. If $N_G(H)=N_{A^{t^{-1}}}(H)\amalg_C N_A(H)$  and $N_A(H)\neq C$ then $[N_A(H):C]=2=[N_{A^{t^{-1}}}(H):C]$ as before (note that either both $N_{A^{t^{-1}}}(H), N_A(H)$ coincide with $C$ or both not, since in every finite quotient they must have the same order).  But if $[N_A(H):C]=2=[N_{A^{t^{-1}}}(H):C]$ then $t$ normalizes $C$ and so this case does not occur.

If on the other hand 
$N_G(H)=HNN(N_{A}(H),C,t)$ and $C\neq N_A(H)$, then factoring out the normal closure of $C$ we get $N_G(H)/\langle\langle C\rangle\rangle= L\amalg \widehat \Z_C$ that contains non-abelian free pro-$\C$ group contadicting solubility of $N_G(H)$. Hence $C=N_A(H)$ in this case.  Then $N_G(H)=C\rtimes \widehat\Z_C$ as required.

\end{proof}

\begin{lemma}\label{pro-C normalizer of cyclic} Let $G$ be a pro-$\C$ group acting on a pro-$\C$ tree $T$ and  $U\cong \widehat\Z_{\C}$ be a procyclic subgroup of $G$ acting freely on $T$. Let $H$ be a non-trivial subgroup of $U$. Then there exists a normal subgroup $K$ of $N_G(H)$ contained in the stabilizer of an edge of $T$ such that $N_G(H)/K$ is  either projective soluble virtually  $\widehat\Z_\C$ acting freely on $T$   or a dihedral pro-$\C$ group $\widehat \Z_C\rtimes C_2$.
\end{lemma}

\begin{proof} Let $N=N_G(H)$ and let $D$ be a minimal $H$-invariant subtree of $T$ (that exists by Proposition \ref{minimal subtree}).  
      Since $D$ is not a vertex, $H$ acts irreducibly on $D$ and so
   by  Proposition \ref{minimal subtree} it is unique. Note that if
   $n\in N$, then $nD$ is also $D$-invariant, and therefore must be equal to
   $D$. 
   Hence $N$ acts irreducibly on $D$ and as $C_N(H)=C_G(H)$, by Lemma 4.2.6(c) in \cite{bib3} (or \cite[Lemma 2.4]{Zal90}) $C_G(H)K/K$ is projective and acts freely on $D$, where $K$ is the kernel of the
   action (and so is the intersection of all stabilizers).

   Since $U\leq C_G(H)$, one deduces from \cite[Proposition 7.4.2]{bib4} that $C_G(H)K/K$ is metacyclic and virtually $\widehat\Z_\C$  (because $UK/K\cong \widehat\Z_\C$) and so $N/K$ is soluble (since $Aut(\widehat\Z_\C)$ is abelian).  Then  by \cite[Proposition 4.2.9]{bib3} (or by \cite[Lemma 2.9]{Zal90}) $N/K$ is either $\Z_\pi\rtimes \Z_\rho$  (where $UK/K\cong \widehat\Z_\C$ is open)  or a dihedral pro-$\C$ group $\widehat \Z_\C\rtimes C_2$.  The case of a profinite Frobenius $N/K$ does not occur since a Frobenius pro-$\C$ group $\Z_\pi\rtimes C_n$ (where $n$ is coprime to $\pi$)   can not contain $\widehat \Z_\C$ as an open subgroup. 
   
   \end{proof}

\begin{proposition}\label{infinite cyclic splitting} Let $G$ be a pro-$\C$ group acting on a pro-$\C$ tree $T$ and  $H\neq 1$ be a cyclic subgroup of $G$ acting freely on $T$. Then  there exists  some normal subgroup $K$ of $N_G(H)$ contained in the stabilizer of an edge such that  $K\leq C_G(H)$,  $C_G(H)/K=P\times H$ is projective and $N_G(H)/K$ is projective-by-abelian. Moreover, if $R$ is a soluble subgroup of $N_G(H)$  then one of the following holds:
	
	\begin{enumerate}
		\item[(i)] $R/K$ is a projective soluble group $\Z_\pi\rtimes \Z_\rho$ ($\pi\cap \rho=\emptyset$)  acting freely on $T$;
		
		\item[(ii)]  $R/K$ is   a profinite Frobenius group $\Z_\pi\rtimes C_n$;
		
		\item[(iii)] $R/K$ is a dihedral pro-$\C$ group $ \Z_\pi\rtimes C_2$.
	\end{enumerate}
	
	In particular one of these items  holds for any vertex stabilizer of $N_G(H)$. Moreover, the image of  any vertex stabilizer of $N_G(H)$ in $N_G(H)/K$ is finite cyclic.
\end{proposition}

\begin{proof} Let $N=N_G(H)$ and let $D$ be a minimal $H$-invariant subtree of $T$ (that exists by Proposition \ref{minimal subtree}).  
	Since $D$ is not a vertex, $H$ acts irreducibly on $D$ and so
	by  Proposition \ref{minimal subtree} it is unique. Note that if
	$n\in N$, then $nD$ is also $D$-invariant, and therefore must be equal to
	$D$.  Hence $N$ acts irreducibly on $D$. 
	
	Let $K$ be the kernel of this
	action. Then $H\cap K=1$ and as both are normal in $N$ we have $KH=K\times H$.
	
	By   \cite[Lemma 4.2.6(c)]{bib3} (or \cite[Lemma 2.4]{Zal90}) ,  $C_G(H)/K$ is projective and acts freely on $D$.  Hence $C_G(H)/K=P\times H$ with $P$ projective coprime to $H$ (indeed as $C_G(H)$ is projective, $C_G(H)/H$ is coprime to $H$ and hence $H$ complements by the profinite version of the Schur-Zassenhaus theorem cf. \cite[Theorem 2.3.15]{bib4}).  
	As $H$ is cyclic,  $N_G(H)/C_G(H)$ is abelian, and therefore so are vertex stabilizers in $N_G(H)/K$ (indeed the image of a vertex stabilizer in $N_G(H)/K$ does not intersect  $C_G(H)/K$  that acts freely on $D$ and so isomorphic to its image in $N_G(H)/C_G(H)$ which is abelian).

	If $R\leq N_G(H)$ is soluble, one deduces from  \cite[Proposition 4.2.9]{bib3} ( or  \cite[Lemma 2.9]{Zal90}) that  $R/K$ is either $\Z_\pi\rtimes \Z_\rho$  or profinite Frobenius $\Z_\pi\rtimes C_n$ or a dihedral pro-$\pi$ group $ \Z_\pi\rtimes C_2$, where $\Z_\pi\rtimes \Z_\rho$, $\Z_\pi$ act freely. In particular, the  vertex stabilizers in $R/K$ are finite cyclic.

\end{proof}

\begin{corollary}\label{not fictitious normalizers} With the hypotheses of Proposition \ref{infinite cyclic splitting} suppose $G=A\amalg_{H} B$ 
	and $N_A(H)\neq N_G(H)\neq N_B(H) $ or $G=HNN(A,H,t)$. Then $N_G(H)$ is metabelian and $N_G(H)=N_A(H)\amalg_C N_B(H)$ with $[N_A(H):H]=2=[N_B(H):H]$ in the first case and  $N_G(H)=H\rtimes \widehat\Z_\C$ with $H=N_A(H)$ in the second case. Moreover,  $N_G(H)/K$ is isomorphic to one of the groups listed in (i)-(iii) of Proposition \ref{infinite cyclic splitting}.
	
\end{corollary}

\begin{proof} By Corollary \ref{centralizer}  $C_G(H)$ contains a free pro-$\C$ group $F$.  As $K$ is cyclic as it is contained in a conjugate of $H$ and $H$ is cyclic),  if $F$ is non-abelian free pro-$\C$, then $C_G(H)/K$ also contains a non-abelian free pro-$\C$ group (since cyclic $K$ must intersect $F$ trivially by \cite[Corollary 8.6.4]{bib4}), contradicting the projectivity of $C_G(H)/K=P\times H$, since $F\times H$ is not projective. Thus $F$ is cyclic and so $N_G(H)$ is soluble. Then $C_G(H)/K$ is soluble and so is a semidirect product of cyclic groups of coprime order (see \cite[Exercise 7.7.8 (b)]{bib4}); such pro-$\C$ group can contain a cyclic free pro-$\C$ group only if it is abelian. It follows that  $C_G(H)/K$ is abelian. The rest follows directly from Corollary \ref{soluble normalizer} and  Proposition \ref{infinite cyclic splitting}.
\end{proof}

\section{Excluding Hyperbolic-Elliptic $\mathbb{Z}_{\widehat{\C}}$-Splitting }\label{sec4}

\bigskip
\begin{proposition}\label{trivial edge stabilizers} Let $G$ be a pro-$\C$ group admitting a cyclic splitting $G=A\amalg_{C} B$ (resp. $G=HNN(A,C, t)$). If $A$ admit an action on an infinite pro-$\C$ tree with trivial edge stabilizers where $C$ is elliptic, then so does $G$.

\end{proposition}

\begin{proof} Suppose  $A$ acts on a pro-$\C$ tree $T$ with trivial edge stabilizers. Let $S$ be the standard pro-$\C$ tree corresponding to the splitting of $G$. Then $G\backslash S$ has one edge only.  Let $\widehat{S}$ be a standard refinement of $S$ with respect to $T$ (by blowing up the vertex whose vertex group is $A$). By  Proposition \ref{refinement} it is a pro-$\C$ tree on which $G$ acts and there exists an edge $e\in E(\widehat{S})$ whose edge stabilizer $G_e=C$ and it is a unique edge with non-trivial stabilizer up to translation. Denote by $v,w$ the vertices of $e$  and let $\Delta=G\{e\cup v\cup w\}$ be the span of translations of $e$. By collapsing connected components of $\Delta$ (see Definition \ref{collapse}) we obtain a pro-$C$ tree  on which $G$ acts with trivial edge stabilizers.   

\end{proof}

\begin{theorem}
   Let $G$ be a  pro-$\C$ group whose every action on  a pro-$\C$ tree with trivial edge stabilizers admits a global fixed point. Then $G$ does not admit a pair of hyperbolic-elliptic cyclic splittings. 
\end{theorem}

\begin{proof}
    We  prove it by contradiction. Suppose the hyperbolic-elliptic case is possible and so there exist two cyclic splittings of $G$  over  cyclic groups $C_1$ and $C_2$, i.e. either $G=A_i\amalg_{C_i} B_i$ or $G=HNN(A_i,C_i, t_i)$ for $i=1, 2$. Let  $T_1$, $T_2$ be the standard pro-$\C$ trees for these splittings and suppose w.l.o.g $c_1$ is hyperbolic in $T_2$ and $c_2$ is elliptic in $T_1$. 

 Since $c_2$ is elliptic in $T_1$, $C_2$ stabilizes a vertex of $T_1$, i.e. the subtree of fixed points (cf. \cite[Theorem 4.1.5]{bib3} or \cite[Theorem 2.8]{ZM}) $T_1^{C_2}\neq \emptyset$. Since $c_1$ is hyperbolic in $T_2$, it intersects trivially any conjugate of $A_2$ and so $A_2$ acts on $T_1$ with trivial edge stabilizers. Then by Proposition \ref{trivial edge stabilizers}   $A_2$ must fix a vertex $v$ in $T_1$.

 Now we will consider  two possible cases.
    
\bigskip
    \textbf{Case 1:} ($G=A_2\amalg_{C_2}B_2$)\\
     By symmetry  $B_2 $ stabilize a vertex in $T_1$ as well. 
  
  Hence $A_2, B_2$ are contained in some conjugate of $A_1$ or $B_1$, say that $A_2\leq A_1$ (note that for the HNN-extension case of the first splitting  there is no $B_1$). Then $B_2$ also has to be in $A_1$ (because otherwsie $A_2\cap B_2=1$ by Proposition \ref{intersection}). Hence $A_2$, $B_2\leq A_1$ and so $G=A_2\amalg_{C_2}B_2\leq A_1$, a contradiction.\\

\textbf{Case 2:} ($G=HNN(A_2,C_2,t_2)$)\\

Since  $A_2$ fixes a vertex in $T_1$,   w.l.o.g  we may assume that $A_2\leq A_1$, and we are going to prove that $t_2$ is in $A_1$. Suppose on the contrary $t_2\notin A_1$. We know that $c_{2}$ and $c^{t_2}_2$ are in $A_{2}<A_{1}$. Then  $c_2\in (A_1)^{t^{-1}_2}$ and so $c_2\in A_1\cap A_1^{t^{-1}_2}$. By Proposition \ref{intersection} $c_2\in A_1\cap A_1^{t^{-1}_2}<C^{a}_1$, for some $a\in A_1$. Then $c_2^{a^{-1}}\in C_1$, which is absurd, because $c_1$ is hyperbolic and so $C_1$ must act freely on $T_2$. 

Thus $A_2$ and $t_2$ are in $A_1$, and $G=HNN(A_2,C_2,t_2)=\langle A_2,t_2\rangle\leq A_1$. So $G\leq A_{1}$, a final contradiction. \\
\end{proof}

Since any pair of abstract cyclic splittings is either elliptic-elliptic or elliptic-hyperbolic or hyperbolic-hyperbolic and by Remark \ref{abstract hyperbolic is profinitely hyperbolic} ellipticity and hypertbolicity of an element survives in the pro-$\C$ completion we deduce the following

\begin{corollary} Let $G_1$ and $G_2$ be residually $\C$ abstract cyclic splittings having the same pro-$\C$ completion $G$  that cannot act on a pro-$\C$ tree with trivial edge stabilizers having no  global fixed points. Then the  pair of corresponding pro-$\C$ splittings of $G$ is either elliptic-elliptic or hyperbolic-hyperbolic.
	
\end{corollary}


\section{Normalizers in the hyperbolic-hyperbolic case }\label{sec5}

\noindent Our objective  in this section is to prove the pro-$\C$ case of Theorem 3.6 of \cite{bib1}. Thus we fix in this section two hyperbolic-hyperbolic splittings  $G=A_i\amalg_{C_i} B_i$ or $G=HNN(A_i,C_i, t_i)$ for $i=1, 2$, and denote by   $T_{1}$ the standard pro-$\C$ tree of the first splitting, and by $T_{2}$ the standard pro-$\C$ tree of the second splitting.\\

\noindent 
Note that since  $c_{1}$ is a hyperbolic element of $T_{2}$ (the standard pro-$\C$ tree of the second splitting) then there exists a unique minimal $C_1$-invariant pro-$\C$ subtree $D_{1}$ of $T_{2}$, such that $C_{1}$ acts irreducibly on $D_{1}$ (see Proposition \ref{minimal subtree}). Hence  $D_{1}$ is $N_{G}(C_{1})$-invariant. Since $C_{1}\triangleleft_{c}N_{G}(C_{1})$ and $C_{1}$ acts irreducibly on $D_{1}$ then $N_{G}(C_{1})$ acts irreducibly on $D_{1}$ (cf. \cite[Remark 4.2.1 (b)]{bib3}).  Denoting by $K_{1}$ the kernel of the action of  $N_{G}(C_{1})$  on $D_{1}$ we deduce that  $N_{G}(C_{1})/K_{1}$ acts irreducibly and faithfully on $D_{1}$. Observe that as $D_{1}\subseteq T_{2}$, then $K_{1}\leq C^{g}_{2}$ for some $g\in G$ (since all edge stabilizers are conjugate of $C_2$, and $K_1$ is contained in each stabilizer of elements from $D_1$) and so w.l.o.g. we may assume that $K_1\leq C_2$. As $K_1$ and $C_1$ are normal in $N_G(C_1)$ and they do not intersect (because $c_1$ is hyperbolic and $K_1$ is elliptic in $T_2$), they commute, i.e. $C_1K_1=C_1\times K_1$. It follows then that $K_1$ normalizes $N_G(C_1)$ and  $C_G(C_1)$. By Proposition \ref{infinite cyclic splitting} $C_G(C_1)/K_1=P\times C_1$ with $P\times C_1$ projective.

Considering the action of $C_2$ on the minimal $C_2$-invariant subtree $D_2$ of $T_1$, and assuming that $K_1\neq 1$ we deduce that $K_1$ acts irreducibly on $D_2$ (because the action is free on $T_1$) and therefore so is $C_1\times K_1$ by \cite[Remark 4.2.1 (b)]{bib3}. Then by \cite[Corollary 4.1.7]{bib3} $C_1$ fixes a vertex in $D_2$ and being normal in $C_1\times K_1$, by \cite[Proposition 4.2.2]{bib3} it is the kernel of the action of it.  Hence $N_G(K_1)$   also acts irreducibly on $D_2$ having $C_1$ as the kernel of the action and as $N_G(C_2)\leq N_G(K_1)$, $N_G(C_2)$ normalizes $C_1$.  The group $C_2$ being normal in $N_G(C_2)$ and intersecting $C_1$ trivially must commute with $C_1$. Thus we proved the following analog of \cite[Lemma 3.3]{bib1}.

\begin{lemma}\label{commute} Let \(G\) be pro-$\C$ group. Consider two hyperbolic-hyperbolic splittings \(G=A_i\amalg_{C_i} B_i\) or \(G=HNN(A_i,C_i, t_i)\) for $i=1, 2$.  If there are $1\neq H_1\le C_1$ and $1\neq H_2\le C_2$  that commute, then $C_1$ and $C_2$ commute.

\end{lemma}

Indeed, if $H_1$ and $H_2$ commute, then the kernel $K_1$ is non-trivial, since $H_2$ must act trivially on $D_1$ by \cite[Corollary 4.1.7 and Proposition 4.2.2]{bib3}.

  Thus as in \cite[Section 3]{bib1}  the study of normalizers splits up into two cases: $C_1, C_2$  commute  up to conjugation or  do not  commute (i.e. $K_1=1$).  

\subsection{$C_1$ and $C_2$ are  commute up to conjugation}

So in this case $K_1\neq 1\neq K_2$ and w.l.o.g. we assume that $C_1$ and $C_2$ commute.

\begin{definition} We will give the following definitions:

\begin{enumerate}
\item[(i)] The semidirect product $C\rtimes \Z/2$, where $C$ is infinite cyclic and $\Z/2$ acts by invertion will be called an infinite dihedral group.

\item[(ii)] A semidirect product $\Z_\pi\rtimes \Z_\rho$, where the action is by inversion  will be called pseudo Klein bottle; it is isomorphic to $(\Z_\pi\times \Z_{\mu})\rtimes \Z_2$.    If $\Z/2\in \C$ the pro-$\C$ completion of the Kleinian bottle is called a pro-$\C$ Klein bottle.
\end{enumerate}

\end{definition}
 
The following Proposition is the pro-$\C$ analog of the Proposition 3.3 of \cite{bib1}.

\begin{proposition} \label{prop1} Let $1\neq H_i\leq C_i$.  Then    $N_{G}(H_{i})$ is soluble and  isomorphic  to one of the following groups:

\begin{enumerate}

\item[(i)] $N_G(H_i)\cong C_i\times \widehat\Z_\C$;
\item[(ii)]  $\Z/2\in \C$ and
$N_{A_i}(H_1)\amalg_{C_i}N_{B_i}$, $[N_{A_i}(H_i):C_i]=2=[N_{B_i}(H_i):C_i]$. Moreover, setting $\pi,\rho$ to be sets of primes such that $C_i\cong \Z_\pi$ and $\rho=\pi\setminus\{2\}$ one of the following holds
\begin{enumerate}
\item[(a)] $N_{A_i}(H_i)$, $N_{B_i}(H_i)$ are torsion free and $N_G(H_i)$ is the pseudo Klein bottle  $(\Z_\rho\times \widehat\Z_\C)\rtimes \Z_2$ with $\Z_2$ acting by inversion, or 
\item[(b)] both $N_{A_i}(H_i)$ and $N_{B_i}(H_i)$ are  infinite dihedral  and
$N_G(H_1)$ is  isomorphic to $(\Z_{\pi}\times \widehat\Z_\C)\rtimes \Z/2$ with $\Z/2$ acting by inversion; 
\item[(c)] One of $N_{A_i}(H_i)$, $N_{B_i}(H_i)$ is torsion free and  the other  is  infinite dihedral, and   $N_G(H_1)\cong (\Z_\rho\times \widehat\Z_\C)\rtimes \Z_2)\rtimes \Z/2)\cong (\Z_\rho\times \widehat\Z_\C)\rtimes (\Z_2\rtimes \Z/2)$.
\end{enumerate}
In particular $N_G(H_i)$ contains an abelian normal subgroup of index $2$.
\end{enumerate}
\end{proposition}

\begin{proof}  
Let $K_i$ be the kernel of the action of $N_G(H_i)$  on $D_i$ which is non-trivial by the hypothesis at the beginning of the subsection and note that $K_i$ acts freely on $T_i$. By Proposition \ref{teo6}  
$N_{G}(H_i)=N_{A_i}(H_i)\amalg_{C_i} N_{B_i}(H_i)$ if $G=A_i\amalg_{C_i}B_i$, and either $N_G(H_i)=N_{A_i^{t_i^{-1}}}(H_i)\amalg_{C_i} N_{A_i}(H_i)$ or $N_G(H_i)=HNN(N_{A_i}(H_i),C_i,t_i)$ if $G=HNN(A_i, C_i, t_i)$.  Moreover, the splittings of the normalizer $N_G(H_i)$ as the amalgated free products can not be fictitious, since otherwise $N_G(H_i)$ stabilisers a vertex in $T_i$ contradicting that $K_i$  acts on $T_i$ freely.  By Corollary \ref{not fictitious normalizers} $N_G(H_i)$ is metabelian and  either $N_G(H_i)\cong C_i\rtimes \widehat\Z_\C$ or 
$N_{A_i}(H_i)\amalg_{C_i}N_{B_i}$ with $[N_{A_i}(C_i):C_i]=[N_{B_i}(C_i):C_i]=2$.

Recall that $[C_1,C_2]=1$ and so $C_2$ is normal in $N_{G}(H_1)$ and $C_1$ is normal in $N_G(H_2)$. If  both $N_G(H_i)\cong C_i\rtimes \widehat\Z_\C$, then  since $C_1\cap C_2=1$, $N_{G}(H_i)$ embeds in $N_{G}(H_i)/C_1\times N_{G}(H_i)/C_2$ and so is abelian. Thus $N_G(H_i)\cong C_i\times \widehat\Z_\C$ and we have case (i).

Otherwise, either $N_G(H_1)$ or $N_G(H_2)$ is generalized dihedral, say $$N_G(H_1)=N_{A_1}(H_1)\amalg_{C_1} N_{B_1}(H_1), [N_{A_1}(C_1):C_1]=[N_{B_1}(C_1):C_1]=2$$ and we have case $(ii)$ for $N_G(H_1)$.  Put $N_1=N_{A_1}(H_1)$, $N_2=N_{B_1}(H_1)$. Let $\pi$ be the set of primes such that $C_1\cong \Z_\pi$ and $\rho=\pi\setminus \{2\}$.

\medskip

If $N_G(H_i)$ is torsion free, then either $N_i\cong \Z_\pi$  or a pseudo Klein bottle group $\Z_{\rho}\rtimes \Z_2$  and so $N_G(H_1)$ is  a pseudo Klein bottle $(\Z_\rho\times \widehat\Z_\C)\rtimes \Z_2$ with $\Z_2$ acting by inversion, i.e.  we have $(ii)$(a).
  
  \smallskip
 If both $N_1$ and $N_2$ are  infinite dihedral  then 
$N_G(H_1)$ is  isomorphic to $(\Z_{\pi}\times \widehat\Z_\C)\rtimes \Z/2$ with $\Z_2$ acting by inversion. Thus 
  $(ii)$(b) holds.
  
  \smallskip
If one of $N_j$, say  $N_2$ is torsion free and  $N_1$ is  infinite dihedral,  then the normal closure $\langle\langle N_2\rangle\rangle_N$ of $N_2$ in $N$ is the pseudo Klein bottle $(\Z_\rho\times \widehat\Z_\C)\rtimes \Z_2$ and so $N_G(H_1)\cong (\Z_\rho\times \widehat\Z_\C)\rtimes \Z_2)\rtimes \Z/2)=(\Z_\rho\times \widehat\Z_\C)\rtimes (\Z_2\rtimes \Z/2)$, so $(ii)$(c) holds. 
  
 Finally we show that $N_1=\Z/2\times C_0$, $C_0\leq C_1$ is impossible. Indeed, in this case considering the action of $N_1$ on $D_1$ we see that $\Z/2$ must be in the kernel, contradicting with $K_1$ being torsion free. 
\end{proof}

\begin{remark}\label{virtually commute}  If the both  cyclic splittings are $\widehat\Z_\C$-splittings, then the subcases (a),(b) and (c) of Proposition \ref{prop1} (ii) are  the pro-$\C$ completion of the Klein bottle, a euclidean 4-branched sphere and an euclidean 2-branched (real) projective plane respectively. To see the last subcase , one observes that  an element $y$ of order $2$ in $N_1$ acts on $C_2$ by inversion (as well as $x$) and so $z=yx$ fixes $C_2$. Thus $N$ has the following pro-$\C$ presentation: $N=\langle z, x\mid zx^2z^{-1}=x^{-2}, xz^2x^{-1}=z^{-2}, (zx^{-1})^2=1\rangle$.   
	
	In particular in all these cases $N$ has normal   subgroup of index 2 isomorphic to $\Z_2\times \Z_2$.
\end{remark}

\begin{corollary} \label{metacyclic p=2} With the hypothesis  of Proposition \ref{prop1} assume  $N_G(H_i)$ is torsion free. Then one of the following holds.
\begin{enumerate}
\item[(ii)]  $N_{G}(H_{i})\cong C_i\times \widehat\Z_\C$.
\item[(iii)] the pseudo Klein bottle $(\Z_\rho\times \widehat\Z_\C)\rtimes \Z_2$ with $\Z_2$ acting by inversion.
\end{enumerate}
\end{corollary}

Assuming  that $2$ is not involved in $\C$ i.e. $\C$ consists of groups of odd order, the description is much nicer and easier. Despite that this case does not  arise from the profinite completion of the abstract  construction., we state it as 

\begin{corollary} \label{2not wnvloved} Suppose $\C$ consists of groups of odd order. Then $N_{G}(H_{i})\cong C_i\times \widehat\Z_\C$.
	
	\end{corollary}

\subsection{$C_1$ and $C_2$ do not virtually commute up to conjugation}

This means that $K_1, K_2$ are trivial.

\begin{proposition} \label{not virtually commute} Let $1\neq H_i\leq C_i$. Then $C_G(H_i)$ is projective and $N_G(H_i)/C_G(H_i)$ is abelian. Moreover, either $N_G(H_i)\subset A_i\cup B_i$  or  $N_{G}(H_{i})$ is soluble and is  isomorphic  to one of the following groups:

\begin{enumerate}
	\item[(i)]  $\Z_\pi\rtimes \Z_\rho$ ($\pi\cap \rho=\emptyset$); 
	\item[(ii)] $\mathbb{Z}_{\pi}\rtimes \Z/2$ infinite dihedral pro-$\C$ group;

	\item[(iii)] a profinite Frobenius group $\Z_\pi\rtimes C_n$

\end{enumerate}
\end{proposition}

\begin{proof} 
	 As $K_i$ is trivial,  by Proposition \ref{infinite cyclic splitting}    $C_G(H_i)=P\times H_i$ is projective and $N_G(H_i)$ is projective-by-abelian. Moreover by Corollary  \ref{not fictitious normalizers} either  $N_G(H_i)\subset A_i\cup B_i$  or  $N_{G}(H_{i})$ is soluble and is  isomorphic  to one of the following groups:
	\begin{itemize}
		\item $\Z_\pi\rtimes \Z_\rho$ ($\pi\cap \rho=\emptyset$)  acting freely on $D_1$;
		\item $ \widehat{\mathbb{Z}}_{\pi}\rtimes \Z/2$ infinite dihedral pro-$\C$ group.
		
		\item $ \Z_\pi\rtimes C_n$ a profinite Frobenius group.
		\end{itemize}

		 Thus we have itens $(i)$, $(ii)$ or $(iii)$.

\end{proof}

\begin{corollary}\label{trivial kernel} 
\begin{enumerate}
\item[(i)] If $N_G(H_1)=HNN(N_{A_1}(H_1),C_1, t)$ (cf. Proposition \ref{teo6}) then $N_G(H_1)\cong\widehat \Z_C$;
\item[(ii)] If all the groups in $\C$ are of odd order then (ii)  does not appear in the statement of Proposition \ref{not virtually commute};
\item[(iii)]
If $C_i\cong \widehat\Z_\C$ then $P=1$ so that $C_G(H_i)=H_i$, $\pi\cup \rho=\pi(C)$ and so $\Z_\pi\rtimes\Z_\rho$ in item(i)  is virtually $\widehat\Z_C$.
\end{enumerate}
\end{corollary}

\begin{corollary} \label{metacyclic} Assume $N_G(H_i)$ is torsion free. Then it is projective and  unless $G=A_i\amalg _{C_i} B_i$ with $N_G(C_i)\subset A_i\cup B_i$, $N_G(C_i)$ is soluble and isomorphic to
  $ \Z_\pi\rtimes \Z_\rho$ ($\pi\cap \rho=\emptyset$).
\end{corollary}

\begin{proof} Put $N=N_G(H_i)$ and consider its action on $D_i$. By Proposition \ref{infinite cyclic splitting} its vertex stabilizer $N_v$ has to be finite cyclic and so is trivial in our case. Thus $N$ acts freely on $D_i$ and therefore is projective (see \cite[Theorem 2.6]{ZM} or \cite[Theorem 4.2.2]{bib3}). Now the result follows from Proposition \ref{not virtually commute}. 
	
	\end{proof}

\section{Hyperbolic-hyperbolic case}\label{sec6}

\begin{proposition}\label{Lemma.pro-2} Let $G$ be a pro-$\C$ group admitting a cyclic splitting $G=A\amalg_{C} B$ (or $G=HNN(A,C,t)$). If $A$ admits an action on a pro-$\C$ tree $T$ with finite cyclic edge stabilizers and no global fixed point, where $C$ is elliptic, then so does $G$.
	
\end{proposition}

\begin{proof}
	 Suppose  $A$ acts on a pro-$\C$ tree $T$ with cyclic edge stabilizers. Let $S$ be the standard pro-$\C$ tree corresponding to the splitting of $G$. Then $G\backslash S$ has one edge only.  Let $\widehat{S}$ be a standard refinement of $S$ with respect to $T$ (see Definition \ref{standard refinement}). By  Proposition \ref{refinement} it is a pro-$\C$ tree on which $G$ acts. If all vertex stabilizers are finite cyclic, we are done. Otherwise,   there exists an edge $e\in E(\widehat{S})$ whose edge stabilizer $G_e=C$ and it is a unique edge that does not have stabilizer of finite order  up to translation. Denote by $v,w$ the vertices of $e$  and let $\Delta=G\{e\cup v\cup w\}$ be the span of translations of $e$. By collapsing connected components of $\Delta$ (see Definition \ref{collapse}) we obtain a pro-$\C$ tree $\widehat{T}$ on which $G$ acts with finite cyclic edge stabilizers.

\end{proof}

\begin{lemma}\label{splitting over normalizer} Let $G=A\amalg_{C} B$ (or $G=HNN(A,C,t)$) be a cyclic splitting. Then $G$ splits as an amalgamated free pro-$\C$  product or as an HNN-extension over $N_{A}(C)$ in one of the following ways:
\begin{itemize}
	\item[(a)] $G=A\amalg_{N_{A}(C)} (N_{G}(C)\amalg_{N_{B}(C)} B)$, if $G=A\amalg_{C} B$;
	\item[(b)] $G=A_1\amalg_{N_{A}(C)}HNN(N_G(C), N_{A^{t^{-1}}}(C),t)$, if $G=HNN(A,C,t)$ and $N_G(C)=N_{A}(C)\amalg_{C} N_{A^{t^{-1}}}(C)$ (cf. Proposition \ref{teo6}) ;
	\item[(c)] $G=A\amalg_{N_{A}(C)}N_{G}(C)$,  if $G=HNN(A,C,t)$ and\\  $N_{G}(C)=HNN(N_{A}(C),C,t)$ (cf. Proposition \ref{teo6}).
\end{itemize} 

\end{lemma}

\begin{proof} Suppose  the splitting is  an amalgamated free pro-$\C$ product $G=A\amalg_{C}B.$ Then by Proposition \ref{teo6} $N_G(C)=N_{A}(C )\amalg_{C} N_{B}(C)$ and therefore $G$ admits a decomposition as follows:

\begin{equation}\label{eq2k}
	\begin{array}{cl} 
	G       &=A\amalg_{C}B=A\amalg_{N_A(C)}(N_A(C)\amalg_CN_B(C))\amalg_{N_B(C)} B\\
	&=A\amalg_{N_{A}(C)} (N_{G}(C)\amalg_{N_{B}(C)} B).\\
	\end{array}
	\end{equation}
Thus $G$ splits as a  free pro-$\C$ product with $N_{A}(C)$ amalgamated  in this case.\\

On the other hand if the first splitting is an HNN-extension  $G=HNN(A,C,t)$, then by Proposition \ref{teo6} the normalizer $N_G(C)$ is an amalgamated free pro-$\C$ product $N_G(C)=N_{A}(C)\amalg_{C} N_{A^{t^{-1}}}(C)$ or  an HNN-extension $N_G(C)=HNN(N_{A}(C),C,t)$. We claim that in the first case 
we have $G=A\amalg_{N_{A}(C)}HNN(N_G(C), N_{A^{t^{-1}}}(C),t)$. Indeed
 we have the decomposition of $G$ as the following graph of groups:

$$
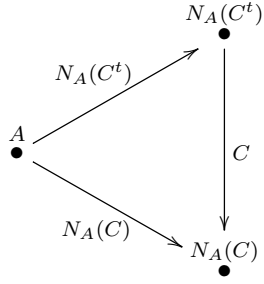
\xymatrix{&&{\overset{N_A(C^t)}{\bullet}}\ar[dd]^C\\
{\overset{A}{\bullet}}\ar[rru]^{N_{A}(C^t)}\ar[rrd]_{N_{A}(C)}&&\\
&&{\overset{N_{A}(C)}{\bullet}}}$$
where the maximal subtree contains two edges on the left.
To see this one can collapse these two fictitious edges on the left to get $HNN(A,C,t)$. On the other hand if we choose lower left and vertical edges to be in the maximal subtree, then the graph of groups is 

$$
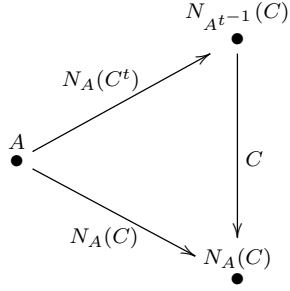
\xymatrix{&&{\overset{N_{A^{t^{-1}}}(C)}{\bullet}}\ar[dd]^C\\
{\overset{A}{\bullet}}\ar[rru]^{N_{A}(C^t)}\ar[rrd]_{N_{A}(C)}&&\\
&&{\overset{N_{A}(C)}{\bullet}}}$$
and so if 
one collapses the vertical edge, one gets

\begin{equation}\label{eq3k}
\begin{array}{cl}
G &=HNN(A\amalg_{N_{A}(C)}N_{G}(C),N_{A}(C^t),t^{-1})\\
&=A\amalg_{N_{A}(C)}HNN(N_G(C), N_{A^{t^{-1}}}(C),t)\\
\end{array}
\end{equation}
where the second equality follows because $N_G(C)=N_{A}(C)\amalg_{C} N_{A^{t^{-1}}}(C)$ contains associated subgroups $N_{A}(C)$ and $N_{A^{t^{-1}}}(C)$ so that considering presentations of both sides one sees they they are equivalent.
   This finishes the proof in this case.

In the second case (recall that in this case $N_G(C)=HNN(N_{A}(C),C,t)$, and so $t$ normalizes $C$)
\begin{equation}\label{eq4k}
	\begin{array}{cl} 
	G       &=HNN(A,C,t)\\
	&=\langle A,t| C^{t}=C\rangle\\
	&=A\amalg_{C} (C\rtimes \langle t\rangle)\\
	&=A\amalg_{N_{A}(C)}N_{A}(C)\amalg_{C}(C\rtimes \langle t\rangle)\\
	&= A\amalg_{N_{A}(C)} HNN(N_{A}(C),C,t)\\
	&=A\amalg_{N_{A}(C)}N_{G}(C)
	\end{array}
	\end{equation}
This finishes the proof.

\end{proof}

\begin{lemma}\label{elliptic factors} Let $G$ be a  pro-$\C$ group  which can not act on a pro-$\C$ tree with finite cyclic edge stabilizers   having no  global fixed point. Let $G=A_{1}\amalg_{C_{1}} B_{1}$ (or $G=HNN(A_{1},C_{1},t_{1})$), and $G=A_{2}\amalg_{C_{2}} B_{2}$ (or $G=HNN(A_{2},C_{2},t_{2})$) be two hyperbolic-hyperbolic cyclic splittings of $G$. Suppose  $[C_1,C_2]=1$. 
\begin{enumerate}
\item[(i)]
If $G$ splits over $N_{A_1}(C_1)$, then $A_2,B_2$ are elliptic in this splitting.
\item[(ii)] If  $G=HNN(A_{1},C_{1},t_{1})$ and $G$ splits over $N_{A_1^{t_{1}^{-1}}}(C_1)$, then $A_2,B_2$ are elliptic in this splitting.
\end{enumerate}
\end{lemma}

\begin{proof} (i) Suppose that $G$ splits over $N_{A_1}(C_1)$ and let  $S_1$ be the  standard pro-$\C$ tree for this splitting. Since  $C_1$ and $C_2$ commute,  $C_2\leq N_G(C_1)$ and so the group $C_2$ is elliptic in $S_1$. As $C_2$ acts freely on $T_1$, $C_2$ does not intersect any conjugate of $A_1$ and so the kernel of the action of $N_{A_1}(C_1)$ on the minimal $C_1$-invariant subtree $D_1$ of $T_2$ is trivial. Therefore  by Proposition \ref{infinite cyclic splitting}  the vertex stabilizers of $N_{A_1}(C_1)$ in $D_1$ are finite cyclic. In particular, $N_{A_1}(C_1)\cap A_2^g$ and so $N_{A_1}(C_1)^g\cap A_2$ are finite cyclic. This means   that   $A_2$-stabilizers of edges of $S_1$ are finite cyclic groups. Therefore by Proposition  \ref{Lemma.pro-2} $G$ would act on a pro-$\C$ tree with finite edge stabilizers, contradicting the hypothesis, unless $A_2$ stabilizes a vertex of $S_1$. Thus   $A_2$ is elliptic in $S_1$. 

Finally if $G=A_{2}\amalg_{C_{2}} B_{2}$, then by symmetry $B_2$ also elliptic in $S_1$.

\medskip
(ii) The proof of (ii) is similar replacing $N_{A_1}(C_1)$ by  $N_{A_1^{t_{1}^{-1}}}(C_1)$ in the argument above.

\end{proof}

\begin{theorem}\label{teo8}
	Let $G$ be a  pro-$\C$ group  which can not act on a pro-$\C$ tree with finite cyclic edge stabilizers   having no  global fixed point.  Let $G=A_{1}\amalg_{C_{1}} B_{1}$ (or $G=HNN(A_{1},C_{1},t_{1})$), and $G=A_{2}\amalg_{C_{2}} B_{2}$ (or $G=HNN(A_{2},C_{2},t_{2})$) be two hyperbolic-hyperbolic cyclic splittings of $G$. Suppose that $C_1$ commutes with a conjugate of $C_2$. Then  $G=N_G(C_1)=N_G(C_2)$ is virtually abelian isomorphic to one of  pro-$\C$ groups listed in (i) or (ii) of Proposition \ref{prop1}. 
	
\end{theorem}

\begin{proof}  If $G=A_1\amalg_{C_{1}}B_1$, then by Proposition \ref{teo6} $N_G(C_1)=N_{A_1}(C_1)\amalg_{C_1} N_{B_1}(C_1)$ and so by Lemma \ref{splitting over normalizer} $G$ admits a decomposition 
		\begin{equation}\label{eq1}
		G=A_{1}\amalg_{N_{A_{1}}(C_{1})} (N_{G}(C_{1})\amalg_{N_{B_{1}}(C_{1})} B_{1})
		\end{equation}

If $G=HNN(A_1,C_1,t_1)$ , then by Proposition \ref{teo6}  one of the following holds

\begin{equation}\label{eq2} 
		G=A_{1}\amalg_{N_{A_{1}}(C_{1})} HNN(N_G(C_1), N_{A_1^{t_{1}^{-1}}}(C_1),t)
		\end{equation}
 \begin{equation}\label{eq3} G=A_1\amalg_{N_{A_1}(C_1)}N_G(C_1)
\end{equation}

 Thus  in all cases $G$ admits an amalgamation decomposition over  $N_{A_1}(C_1)$.		

\bigskip		
Let $S_1$ be the  standard pro-$\C$ tree for one of these decompositions. Since $C_1$ commute with a conjugate of $C_2$, we may assume w.l.o.g. that $C_1$ and $C_2$ commute. Then by Lemma \ref{elliptic factors}(i)  $A_2$ is contained up to conjugation either in $A_1$ or in the second factor $R$ of one of the decompositions \eqref{eq1}-\eqref{eq3}, so that w.l.o.g we may assume that $A_2\leq A_1$ or $A_2\leq  R$ and similarly $B_2\leq A_1$ or $B_2\leq  R$.

\smallskip		 
We claim that \eqref{eq1}-\eqref{eq3} are fictitious. Suppose not.		
		 If the second splitting is $G=A_{2}\amalg_{C_{2}} B_{2}$, then by Proposition \ref{intersection} $C_2=A_2\cap B_2\leq (N_{A_{1}}(C_1))^{h}$ for some  $h\in G$ contradicting the hypothesis that $c_2$ is hyperbolic in $T_1$.  If the second splitting is an HNN-extension, then by Proposition \ref{intersection} $C_2=A_2\cap A_2^{t_{2}^{-1}}\leq (N_{A_{1}}(C_1))^h$ or $C_2\leq  (N_{A_{1}^{t_2^{-1}}}(C_1))^{h}$  for some  $h\in G$, contradicting  the hyperbolicity of $c_2$ in $T_1$ again.
		  This proves that $A_1= N_{A_1}(C_1)$.
		  
		 \medskip 
		  Thus  we deduce that $G=N_G(C_1)$ if decomposition  \eqref{eq3} holds.
	
	\medskip

Recall that \eqref{eq1} holds, if the first splitting is $G=A_1\amalg_{C_1} B_1$. Then  we not only have   $A_1=N_{A_1}(C_1)$ but  swapping $A_1$ and $B_1$ in the first splitting of $G$  we also have  $N_{B_{1}}(C_{1})=B_1$. Thus $G=N_{G}(C_1)=N_{A_1}(C_1)\amalg_{C_1}N_{B_1}(C_1)$. 

		\bigskip
		
		 Suppose now the decomposition \eqref{eq2} holds and so $G=HNN(N_G(C_1), N_{A_1^{t_{1}^{-1}}}(C_1),t)$.
		 Note that  $[N_{A_1}(C_1) : C_1]\leq 2\geq [N_{A_1^{t^{-1}}}(C_1):C_1]$ since  $N_G(C_1)=N_{A_1}(C_1)\amalg_{C_1} N_{A_1^{t_{1}^{-1}}}(C_1)$ is virtually
		  abelian by Proposition \ref{prop1} and otherwise by Corollary \ref{centralizer} it will be not. Then denoting by $S'_1$  the  standard pro-$\C$ tree for  the decomposition $G=HNN(N_G(C_1), N_{A_1^{t_{1}^{-1}}}(C_1),t)$ we deduce from Lemma \ref{elliptic factors}(ii) that $A_2$ fixes a vertex of $S'$ and hence  is conjugate into $N_G(C_1)$. Thus w.l.o.g we may assume that $A_2\leq N_G(C_1)$. 
		  
		  If the second splitting is an amalgamation then by Proposition \ref{intersection} $C_2=A_2\cap B_2\leq ((N_{A_{1}^{t_{1}^{-1}}}(C_1))^{h}$ for some  $h\in G$,  contradicting the hyperbolicity of $c_2$ in $T_1$.

		 So the second splitting is an HNN-extension and  by Proposition \ref{intersection} $C_2=A_2\cap A_2^{t_{2}^{-1}}\leq ((N_{A_{1}^{t_{1}^{-1}}}(C_1)\cup N_{A_{1}^{t_{1}^{-1}t_{2}^{-1}}}(C_1))^{h}$  
		  for some  $h\in G$ contradicting the hyperbolicity of $c_2$ in $T_1$ again. Therefore this case does not occur.
		  
		   Thus, $G=N_G(C_1)$ in all the cases and  we have either Case $(i)$ or  Case $(ii)$  of Proposition \ref{prop1}.
		
\end{proof}

\bigskip
If the  pro-$\C$ splittings are $\widehat\Z_\C$-splittings then we can apply Corollary \ref{not fictitious normalizers} instead of Proposition \ref{infinite cyclic splitting} to deduce the following  

\begin{corollary}\label{over pro-C}
	Let $G$ be a  pro-$\C$ group  which can not act on a pro-$\C$ tree with finite cyclic edge stabilizers   having no  global fixed point.  Let $G=A_{1}\amalg_{C_{1}} B_{1}$ (or $G=HNN(A_{1},C_{1},t_{1})$), and $G=A_{2}\amalg_{C_{2}} B_{2}$ (or $G=HNN(A_{2},C_{2},t_{2})$) be two hyperbolic-hyperbolic $\widehat\Z_{\C}$-splittings of $G$. Suppose $C_1$ commutes with  $C_2$ up to conjugation. Then  $G=N_G(C_1)=N_G(C_2)$ is virtually abelian isomorphic to one of  pro-$\C$ groups listed in Remark \ref{virtually commute}. 
	
\end{corollary}

If $\C$ consists of groups of odd order Theorem 5.2 simplifies. We state it as

\begin{corollary}\label{odd oder} Suppose $\C$ consists of groups of odd order and let $G$ be a  pro-$\C$ group  which can not act on a pro-$\C$ tree with finite cyclic edge stabilizers   having no  global fixed point.  Let $G=A_{1}\amalg_{C_{1}} B_{1}$ (or $G=HNN(A_{1},C_{1},t_{1})$), and $G=A_{2}\amalg_{C_{2}} B_{2}$ (or $G=HNN(A_{2},C_{2},t_{2})$) be two hyperbolic-hyperbolic cyclic splittings of $G$. Suppose that $C_1$ commutes with a conjugate of $C_2$ . Then $G=N_G(C_1)=N_G(C_2)\cong C_i\times \widehat\Z_\C$. \end{corollary}

\section{2-acylindricity}

In this section we show that the hyperbolic-hyperbolic splittings from Section 4.2 (i.e. over cyclic groups that do not virtually commute) act 2-acylindrically on it standard  pro-$\C$ tree, i.e.  amalgamated subgroups (resp. one of associated subgroups) are malnormal in one of the factors (resp. in the base group). 

\medskip
We shall start with simple lemma. 

\begin{lemma}\label{normalizer of cyclic} Let $G$ be a profinite group and $C$ a cyclic subgroup of $G$. If $1\neq H=C^g \cap C$ for some $g\in G$, then $g\in N_G(H)$.
\end{lemma}

\begin{proof} If $g\not\in N_G(H)$ then there exists a finite quotient of $G$, where the image of $g$ does not belong to the image of $N_G(H)$, so it suffices to prove the statement for finite $G$. But then $H^g$ has the same order as $H$ and both are subgroups of finite cyclic group $C^g$, so $H=H^g$ as required.

\end{proof}

\begin{theorem}\label{teo8}
	Let $G$ be a  pro-$\C$ group  which can not act on a pro-$\C$ tree with finite cyclic edge stabilizers   having no  global fixed point.  Let $G=A_{1}\amalg_{C_{1}} B_{1}$ (or $G=HNN(A_{1},C_{1},t)$), and $G=A_{2}\amalg_{C_{2}} B_{2}$ (or $G=HNN(A_{2},C_{2},t_{2})$) be two hyperbolic-hyperbolic cyclic splittings of $G$ (with $C_1,C_2$ non-trivial). Suppose that $C_1$ does not commutes with a conjugate of $C_2$. Then  $C_1$ is malnormal in $A_1$ or $B_1$ (resp. $C_1$ or $C_2$ is malnormal in  $A_1$). 
	
\end{theorem}

 \begin{proof}

 	Case 1. $G=A_1\amalg_{C_1}B_1$. By contradiction assume that $C_1$  is not malnormal in $A_1$ and $B_1$. Then there exist $a\in A_1-C_1$, $b\in B_1-C_1$ such that $C_1^{a}\cap C_1\neq 1\neq C_1\cap C_1^{b}$.  Then by Lemma \ref{normalizer of cyclic} $H=C_1\cap C_1^a\cap C_1^b$ is normalized by $a$ and $b$ and so $N_G(H)=N_{A_1}(H)\amalg_{C_1} N_{B_1}(H)$ (see Proposition \ref{teo6}) is a non-fictitious splitting, i.e.  $
 	N_{A_1}(H)\neq N_G(H)\neq N_{B_1}(H)$).
 	 Therefore by Corollary \ref{not fictitious normalizers} $N_G(H)$ is metabelian.  Then by Corollary \ref{soluble normalizer}  $[N_{A_1}(H):C_1]=2=[N_{B_1}(H):C_1] $. It follows that $C_1$ is normal, $N_G(C_1)/C_1\cong \Z/2\amalg \Z/2\cong \widehat\Z_C\rtimes \Z/2$ and so $C_G(C_1)$ contains a subgroup isomorphic to $C_1\times \widehat\Z_\C$ (since the image of $\Z/2\amalg \Z/2$ in $Aut(C_1)$ is finite). But this contradicts to projectivity of $C_G(C_1)$ established in  Proposition \ref{infinite cyclic splitting}. So $C_1$ must be trivial,  contradicting the hypothesis. \\
 	
 	Case 2.	
 	Now suppose that $G=HNN(A_1,C_1,t)$. Assume on the contrary that  $C_1$ and $C_1^{t}$ are both not  malnormal  in $A_1$. Then there exist $a_1\in A_1-C_1$, $a_2\in A_1-C_1^{t}$  such that $C_1^{a_1}\cap C_1\neq 1\neq C_1^{ta_2}\cap C_1^{t} $ and so $C_1^{a_1}\cap C_1 \neq 1\neq C_1^{ta_2t^{-1}}\cap C_1$. By Lemma \ref{normalizer of cyclic},  $H= C_1\cap C_1^{a_1}\cap  C_1^{ta_2t^{-1}} $ normalized by $a_1$and $ta_2t^{-1}$. 
 	By Proposition \ref{teo6} $N_{G}(H)$ is either a free amalgamated product $$N_{G}(H)=N_{A_1}(H)\amalg_{C_1}N_{A_1^{t^{-1}}}(H)$$ or an HNN-extension $$N_{G}(H)=HNN(N_{A_1}(H),{C_1},t).$$ 
 	In the first case observing that $a_2\in N_{A_1^{t^{-1}}}(H)$ we see that  $$N_{G}(H)=N_{A_1}(H)\amalg_{C_1}N_{A_1^{t^{-1}}}(H)$$   is  non-fictitious, i.e. $
 	N_{A_1}(H)\neq N_G(H)\neq N_{A_1^{t^{-1}}}(H)$.  Therefore by Corollary \ref{not fictitious normalizers} $N
 	_G(H)$ is soluble.  Then by Corollary \ref{soluble normalizer}  $[N_A(H):C_1]=2=[N_{A^{t^{-1}}}(C_1):C_1] $. Then as in the first case, it follows that $C_1$ is normal, $N_G(C_1)/C_1\cong \Z/2\amalg \Z/2\cong \widehat\Z_C\rtimes \Z/2$ and so $C_G(C_1)$ contains a subgroup isomorphic to $C_1\times \widehat\Z_\C$ (since the image of $\Z/2\amalg \Z/2$ in $Aut(C_1)$ is finite). But this contradicts to projectivity of $C_G(C_1)$ established in  Proposition \ref{infinite cyclic splitting}. So $C_1$ must be trivial,  contradicting the hypothesis.
 	
 \medskip	
 	In the second case $$N_{G}(H)=HNN(N_{A_1}(H),C_1,t)\cong C_1\rtimes \widehat\Z_C$$ is soluble by Corollary \ref{not fictitious normalizers}. This however contradicts Corollary \ref{trivial kernel}(i), since  $C_1\neq 1$.

 \end{proof}

\end{document}